\documentclass[12pt]{article}

\usepackage{graphicx}
\usepackage{amsmath,amssymb,amsfonts,amsthm}
\usepackage{color}
\usepackage[english]{babel}
\usepackage[latin1]{inputenc}
\usepackage{times}
\usepackage[T1]{fontenc}

\newtheorem{thm}{\bf{Theorem}}[section]
\newtheorem{lem}[thm]{\bf{Lemma}}

\newtheorem{df}[thm]{\bf{Definition}}
\newtheorem{cor}[thm]{\bf{Corollary}}
\newtheorem{rem}[thm]{\bf{Remark}}

\newtheorem{prop}[thm]{\bf{Proposition}}

\numberwithin{equation}{section}

\newcommand{\dom}{\operatorname{dom}}
\newcommand{\gph}{\operatorname{gph}}

\newcommand{\Lip}{\operatorname{Lip}}

\newcommand{\argmin}{\operatornamewithlimits{argmin}}

\begin{document}
\title{The NC-Proximal Average for Multiple Functions \thanks{Research by these authors was supported by UBC UGF and by NSERC of Canada.}}
\author{W. Hare\thanks{Mathematics, University of British Columbia, Kelowna. warren.hare@ubc.ca} \and C. Planiden\thanks{Mathematics, University of British Columbia, Kelowna. chayneplaniden@hotmail.com}}
\date{January 15, 2014}
\maketitle
\begin{abstract}
The NC-proximal average is a parametrized function used to continuously transform one proper, lsc, prox-bounded function into another. Until now it has been defined for two functions. The purpose of this article is to redefine it so that any finite number of functions may be used. The layout generally follows that of \cite{proxave}, extending those results to the more general case and in some instances giving alternate proofs by using techniques developed after the publication of that paper. We conclude with an example examining the discontinuity of the minimizers of the NC-proximal average.
\end{abstract}


\section{Introduction}\label{intro}
In 2008, Bauschke, Lucet, and Trienis, first addressed the question of how to transform one convex function into another in a continuous manner \cite{howto}. Given proper convex functions $f_0$ and $f_1$, their proposed solution, the \emph{proximal average}, used Fenchel conjugates to define a parameterized function $PA(x,\lambda)$ such that $PA$ is epi-continuous with respect to $\lambda,$ and $PA(x,0)=f_0(x)$, $PA(x,1)=f_1(x)$ for all $x.$ The proximal average has been studied extensively since its original conception, and many favourable properties and applications of this approach have arisen \cite{proxbas,primdual,respos,compave,autoconj,convhull,optval,proxsadd,proxave,convprox,resave,whatshape,selfreg}. For example, the minimizers of the proximal average function change continuously with respect to $\lambda$ \cite{optval}.

The proximal average has also been generalized and reformulated in a number of useful manners.  For example, in \cite{proxbas}, the proximal average is generalized to a finite number of convex functions.  In \cite{kernave}, the proximal average is generalized to allow for alternate kernels, which further allowed for applications with monotone operators. In \cite{proxsadd}, the proximal average is reformulated to apply to saddle functions.   And, in \cite{proxave}, the proximal average was reformulated to work with two (nonconvex), proper, lsc, prox-bounded functions.   This document generalizes the work done in \cite{proxave} to allow for a finite number of such functions.\smallskip\\
Given two proper, lsc, prox-bounded functions, $f_0$ and $f_1$, the NC-proximal average was originally defined as
    $$PA_r(x,\lambda):=-e_{r+\lambda(1-\lambda)}\left(-(1-\lambda)e_rf_0-\lambda e_rf_1\right)(x)$$
where $\lambda\in[0,1]$ and $e_rf$ is the \emph{Moreau envelope} of $f$ using the \emph{prox-parameter} $r,$ defined as
    $$e_rf(x):=\inf\limits_y\left\{f(y)+\frac{r}{2}|y-x|^2\right\}.$$
Associated with the Moreau envelope, and closely related to the NC-proximal average, is the \emph{proximal point mapping} $P_r f$ defined as
    $$P_rf(x):=\argmin\limits_y\left\{f(y)+\frac{r}{2}|y-x|^2\right\}.$$
In \cite{proxave} the function $PA_r$ is analyzed and a number of propositions and theorems are developed in order to describe its properties. Here, we extend those results for a finite number of proper, lsc, prox-bounded functions $f_i$, $i\in\{1,2,\ldots,m\}.$ We begin by defining the NC-proximal average as
    \begin{equation}\label{eq:DEF} PA_{r,\delta}(x,\lambda):=-e_{r+\delta(\lambda)}\left(-\sum\limits_{i=1}^m\lambda_ie_rf_i\right)(x),\end{equation}
    $$\lambda\in\Lambda:=\left\{(\lambda_1,\lambda_2,\ldots,\lambda_m)\in\mathbb{R}^m:\lambda_i\geq0\mbox{ for all }i\mbox{ and }\sum\limits_{i=1}^m\lambda_i=1\right\},$$
and $\delta$ is any continuous function such that $\delta(\lambda)=0$ if $\lambda=e_i$ (the canonical unit vector whose $i^{th}$ component is 1) for some $i,$ and $\delta(\lambda)>0$ otherwise. This definition generalizes that of \cite{proxave} in two respects. First, the original definition is restricted to outer prox-parameter $r+\lambda(1-\lambda),$ when in fact the $\lambda(1-\lambda)$ term can be replaced by any function $\delta$ as described above. Second, the results found in \cite{proxave} are reworked in order to accommodate any finite number of functions.

\begin{rem}It should be clear that the choice of the function $\delta$ used in defining the NC-proximal average will have a great impact on the parameterized function $PA_{r,\delta}$.  However, it will become clear in this paper that the underlying properties of $PA_{r,\delta}$ are in fact not effected by $\delta$.  As such, for ease of notation, except when necessary we shall simplify $PA_{r,\delta}$ to $PA_r$.
\end{rem}

The remainder of this article is organized as follows. Section \ref{prelim} provides definitions and shows that $PA_r$ is well-defined. Section \ref{prox} explores the prox-regularity and para-prox-regularity aspects of the function, and Section \ref{stab} considers its stability. We conclude, in Section \ref{example}, with some discussion on the minimizers of the NC-proximal average, including an example that demonstrates that the minimizers of the NC-proximal average may be multi-valued and discontinuous.

\section{Preliminaries}
\label{prelim}
Throughout this paper, we use $q$ to represent the norm-squared function, $q(x)=|x|^2$.  This section restates some definitions we need, and shows that under basic assumptions, $PA_r$ is a well-defined function.
\begin{df}
A proper function $f:\mathbb{R}^n\rightarrow\mathbb{R}\cup\{\infty\}$ is said to be \emph{prox-bounded} if there exist $r>0$ and a point $\bar{x}$ such that $e_rf(\bar{x})>-\infty.$ The infimum of the set of all such $r$ is called the \emph{threshold of prox-boundedness}.
\end{df}
\begin{df}
A function is \emph{lower-$\mathcal{C}^2$} on an open set $V$ if it is finite-valued on $V$ and at any point $x\in V$ the function appended with a quadratic term is convex on some open convex neighborhood $V'$ of $x.$ The function is said to be lower-$\mathcal{C}^2$ (with no mention of $V$) if $V=\mathbb{R}^n.$
\end{df}
Our first task is to confirm that $PA_r$ is a well-defined and well-behaved function. The following proposition generalizes \cite[Prop 2.5]{proxave}.
\begin{prop}
\label{prop1}
For $i\in\{1,2,\ldots,m\}$ let $f_i:\mathbb{R}^n\rightarrow\mathbb{R}\cup\{\infty\}$ be proper, lsc, prox-bounded functions with respective thresholds $\bar{r}_i$. Let $r>\max\limits_i\{\bar{r}_i\}.$ Then for all $\lambda\in\Lambda$, $PA_r$ is a proper function in $x$. Furthermore, if $\lambda_i\neq1$ for all $i$, then $PA_r$ defines a lower-$\mathcal{C}^2$ function in $x$. Finally, if for some $i$ one has that $f_i+\frac{r}{2}q$ is convex, then $PA_r(\cdot,e_i)=f_i.$
\end{prop}
\textbf{Proof:} We know that $-e_rf_i$ is well-defined for all $i$, since $r>\bar{r}_i$ for all $i$. By \cite[Lem 2.4]{proxave}, which is extendible to the case of $m$ functions, we know that $-\sum\limits_{i=1}^m\lambda_i e_r f_i$ is a proper, lower-$\mathcal{C}^2$, prox-bounded function, with threshold $\bar{r}\leq\sum\limits_{i=1}^m\lambda_ir=r$. Thus the Moreau envelope of $-\sum\limits_{i=1}^m\lambda_ie_rf_i$ is well-defined and proper whenever the prox-parameter is greater than or equal to $r$ (as is the case when $\lambda\in\Lambda$), and it is lower-$\mathcal{C}^2$ whenever the prox-parameter is strictly greater than $r$ (as is the case when $\lambda\in\Lambda$ and $\lambda_i\neq1$ for all $i$). The last statement is proved by applying \cite[Ex 11.26 (d)]{rockwets} to $PA_r(x,e_i)=-e_r(-e_rf_i)(x)$.\qed

\section{Prox-Regularity}
\label{prox}
In this section, we wish to establish the conditions under which the function $\sum\limits_{i=1}^m\lambda_i e_r f_i$ is para-prox-regular, so that in Section \ref{stab} we may explore the stability of $PA_r.$ Let us recall what we mean by prox-regularity and para-prox-regularity of a function.
\begin{df}
A proper function $f$ is \emph{prox-regular} at a point $\bar{x}$ for $\bar{v}\in\partial f(\bar{x})$ if $f$ is locally lsc at $\bar{x}$ and there exist $\epsilon>0$ and $r>0$ such that
\begin{equation}
\label{eq:prox}
f(x')\geq f(x)+\langle v,x'-x\rangle-\frac{r}{2}|x'-x|^2
\end{equation}
whenever $x'\neq x,$ $|x'-\bar{x}|<\epsilon,$ $|x-\bar{x}|<\epsilon,$ $|f(x)-f(\bar{x})|<\epsilon,$ $v\in\partial f(x),$ and $|v-\bar{v}|<\epsilon.$  We say the function is \emph{continuously prox-regular} at $\bar{x}$ for $\bar{v}$ if, in addition, $f$ is continuous as a function of $(x,v)\in\gph\partial f$ at $(\bar{x},\bar{v}).$ The function is said to be prox-regular at $\bar{x}$ (with no mention of $\bar{v}$) if it is prox-regular at $\bar{x}$ for every $\bar{v}\in\partial f(\bar{x}),$ and simply prox-regular (with no mention of $\bar{x}$) if it is prox-regular at $\bar{x}$ for every $\bar{x}\in\dom f.$
\end{df}
From a graphical point of view, a prox-regular function is one that is locally bounded below by quadratics of equal curvature. Para-prox-regularity is an extension of this idea that includes an extra parameter $\lambda$.
\begin{df}
A proper, lsc function $f:\mathbb{R}^n\times\mathbb{R}^s\rightarrow\mathbb{R}\cup\{\infty\}$ is \emph{parametrically prox-regular} in $x$ at $\bar{x}$ for $\bar{v}\in\partial_xf(\bar{x},\bar{\lambda})$ with compatible parametrization by $\lambda$ at $\bar{\lambda}$ (also refered to as \emph{para-prox-regular} in $x$ at $(\bar{x},\bar{\lambda})$ for $\bar{v}$), with parameters $\epsilon>0$ and $r>0,$ if
\begin{equation}
\label{eq:para}
f(x',\lambda)\geq f(x,\lambda)+\langle v,x'-x\rangle-\frac{r}{2}|x'-x|^2
\end{equation}
whenever $x'\neq x,$ $|x'-\bar{x}|<\epsilon,$ $|x-\bar{x}|<\epsilon,$ $|f(x,\lambda)-f(\bar{x},\bar{\lambda})|<\epsilon,$ $|\lambda-\bar{\lambda}|<\epsilon,$ $v\in\partial_xf(x,\lambda),$ and $|v-\bar{v}|<\epsilon.$
It is \emph{continuously para-prox-regular} in $x$ at $(\bar{x},\bar{\lambda})$ for $\bar{v}$ if, in addition, $f$ is continuous as a function of $(x,\lambda,v)\in\gph\partial_xf$ at $(\bar{x},\bar{\lambda},\bar{v}).$
If the parameter $\bar{\lambda},$ the subgradient $\bar{v},$ or the point $\bar{x}$ is omitted, then the para-prox-regularity of $f$ is understood to mean for all $\bar{\lambda}\in\dom f(\bar{x},\cdot),$ for all $\bar{v}\in\partial_xf(\bar{x},\bar{\lambda}),$ or for all $\bar{x}\in\dom f(\cdot,\bar{\lambda}),$ respectively.
\end{df}
\begin{prop}
\label{cor1}
For $i\in\{1,2,\ldots,m\}$ let $f_i:\mathbb{R}^n\rightarrow\mathbb{R}\cup\{\infty\}$ be proper, lsc, and prox-bounded with threshold $r_i.$ Let $r>r_i$ for all $i.$ Define
$$F(x,\lambda)=\begin{cases}
-\sum\limits_{i=1}^m\lambda_ie_rf_i(x) & ,~\lambda\in\Lambda\\
\infty & ,~\lambda\not\in\Lambda.\end{cases}$$
Then $F$ is continuously para-prox-regular at any $\bar{x}$, with compatible parametrization by $\lambda$ at any $\bar{\lambda}\in\Lambda.$ Moreover, $F$ is lower-$\mathcal{C}^2$ and strictly continuous, and if $(0,y)\in\partial^\infty F(\bar{x},\bar{\lambda})$ then $y=0.$
\end{prop}
\textbf{Proof:} Since $f_i$ is proper, lsc and prox-bounded for all $i,$ \cite[Ex 10.32]{rockwets} gives us that $-e_rf_i$ is lower-$\mathcal{C}^2$ for all $i.$ The sum of lower-$\mathcal{C}^2$ functions is lower-$\mathcal{C}^2,$ and any lower-$\mathcal{C}^2$ function is strictly continuous \cite[Thm 10.31]{rockwets}, so $F$ is lower-$\mathcal{C}^2$ and strictly continuous. Finally, \cite[Thm 9.31]{rockwets} states that strict continuity of $F$ at $(\bar{x},\bar{\lambda})$ is equivalent to $\partial^\infty F(\bar{x},\bar{\lambda})=\{0\},$ which gives us that $(0,y)\in\partial^\infty F(\bar{x},\bar{\lambda})\Rightarrow y=0.$ This gives us all the conditions of \cite[Thm 5.7]{parapr}, and its conclusion is the result we seek.\qed
\begin{rem}
The proof of \cite[Lemma 3.3]{proxave} can also be adapted for a longer, but more direct proof of Proposition \ref{cor1}.
\end{rem}

\section{Stability}
\label{stab}
We are now ready to explore the stability of the NC-proximal average. By Proposition \ref{cor1}, we can see that $PA_r$ is the Moreau envelope of a para-prox-regular function. This allows us to take advantage of the work done in \cite{proxmap}, where the tilt stability and full stability of Moreau envelopes and proximal mappings of para-prox-regular functions was studied.
\begin{thm}\cite[Thm 4.6]{proxmap}
\label{thm1}
Let $F:\mathbb{R}^n\times\mathbb{R}^s\rightarrow\mathbb{R}\cup\{\infty\}$ be proper, lsc, and continuously para-prox-regular at $(\bar{x},\bar{\lambda})$ for $\bar{v}\in\partial_xF(\bar{x},\bar{\lambda}),$ with parameters $\epsilon$ and $r.$ Assume further that $F$ is prox-bounded with threshold $\rho,$ and that $F$ satisfies the following:
\begin{enumerate}
\item $(0,y)\in\partial^\infty F(\bar{x},\bar{\lambda})\Rightarrow y=0,$
\item $(0,\lambda')\in D^*(\partial_xF)(\bar{x},\bar{\lambda}|\bar{v})(0)\Rightarrow\lambda'=0,$
\item $(x',\lambda')\in D^*(\partial_xF)(\bar{x},\bar{\lambda}|\bar{v})(v'), v'\neq0\Rightarrow\langle x',v'\rangle>-\rho'|v'|^2$ for some $\rho'>0,$
\item $\partial_xF(\bar{x},\cdot)$ has a continuous selection $g$ near $\bar{\lambda},$ with $g(\bar{\lambda})=\bar{v}.$
\end{enumerate}
If $\bar{r} > \max\{\rho, \rho', r\}$,
then there exist $K>0$ and a neighborhood $\mathcal{B}=B_\delta(\bar{x}+\frac{\bar{v}}{r},\bar{\lambda},\bar{r})$ such that for all $(x,\lambda,r),(x',\lambda',r')\in\mathcal{B}$ we have that $P_rF_\lambda(x)$ and $P_{r'}F_{\lambda'}(x')$ are single-valued, with
$$|P_rF_\lambda(x)-P_{r'}F_{\lambda'}(x')|\leq K|r(x-\bar{x})-r'(x'-\bar{x}),\lambda-\lambda',r-r')|,$$
where $F_\lambda(x)=F(x,\lambda).$
\end{thm}
\begin{lem}\cite[Lem 4.4]{proxave}
\label{lem4}
Suppose the function $H:\mathbb{R}^n\times\mathbb{R}^s\rightarrow\mathbb{R}\cup\{\infty\}$ is finite, single-valued, and Lipschitz continuous in $(x,\lambda)$ near $(\bar{x},\bar{\lambda})$ with local Lipschitz constant $\Lip H.$ Then
$$(0,\lambda')\in D^*H(\bar{x},\bar{\lambda}|H(\bar{x},\bar{\lambda}))(0)\Rightarrow\lambda'=0,$$
and for $\rho>\Lip H$ one has
$$(x',\lambda')\in D^*H(\bar{x},\bar{\lambda}|H(\bar{x},\bar{\lambda}))(v'),v'\neq0\Rightarrow\langle x',v'\rangle>-\rho|v'|^2.$$
\end{lem}
The next proposition is an analog of \cite[Prop 4.5]{proxave}, rewritten to work with a finite number of functions. The proof of \cite[Prop 4.5]{proxave} is easily adaptable to this setting, so we present only the key details.
\begin{prop}
\label{prop4}
For $i\in\{1,2,\ldots,m\}$, let $f_i:\mathbb{R}^n\rightarrow\mathbb{R}\cup\{\infty\}$ be proper, lsc, and prox-bounded with threshold $r_i.$ Let $r>\max\limits_i\{r_i\},$ and define
$$F(x,\lambda):=-\sum\limits_{i=1}^m\lambda_ie_rf_i(x).$$
If $P_rf_i$ is single-valued and Lipschitz continuous for all $i,$ then the following three properties hold:
\begin{enumerate}
\item $(0,\lambda')\in D^*(\partial_xF)(\bar{x},\bar{\lambda}|\bar{v})(0)\Rightarrow\lambda'=0,$
\item for some $\rho>0$ we have $(x',\lambda')\in D^*(\partial_xF(\bar{x},\bar{\lambda}|\bar{v})(v'),v'\neq0\Rightarrow\langle x',v'\rangle>-\rho|v'|^2,$ and
\item the set-valued mapping $\partial_xF(\bar{x},\cdot)$ has a continuous selection $g$ near $\bar{\lambda}.$
\end{enumerate}
\end{prop}
\textbf{Proof:} Since $P_rf_i$ is Lipschitz continuous, we have that $e_rf_i\in\mathcal{C}^{1+}$ with $\nabla e_rf_i=r(I-P_rf_i)$ \cite[Thm 2.4]{proxmap}. Hence,
\begin{align*}
\partial_xF(\bar{x},\lambda) & =\nabla_x(-\sum\limits_{i=1}^m\lambda_ie_rf_i)(\bar{x},\lambda)\\
 & =r\left[\left(\sum\limits_{i=1}^m\lambda_iP_rf_i(\bar{x})\right)-\bar{x}\right]
\end{align*}
which is linear in $\lambda,$ showing Property 3. Since $P_rf_i$ is single-valued and Lipschitz continuous, we have $\partial_xF(x,\lambda)$ single-valued and Lipschitz continuous. Properties 1 and 2 follow by applying Lemma \ref{lem4}.\qed
\begin{prop}
\label{prop5}
For $i\in\{1,2,\ldots,m\}$, let $f_i:\mathbb{R}^n\rightarrow\mathbb{R}\cup\{\infty\}$ be proper, lsc, and prox-bounded with threshold $r_i.$ Let $r>\max\limits_i\{r_i\}.$ Then $PA_r(\cdot,\lambda)+\frac{r+\delta(\lambda)}{2}q(\cdot-\bar{x})$ is convex for any $\bar{x}.$ Hence, $PA_r(\cdot,\lambda)$ is lower-$\mathcal{C}^2.$
\end{prop}
\textbf{Proof:} Define $F_\lambda:=-\sum\limits_{i=1}^m\lambda_ie_rf_i.$ Then
$$PA_r+\frac{r+\delta(\lambda)}{2}q=-e_{r+\delta(\lambda)}(F_\lambda)+\frac{r+\delta(\lambda)}{2}q.$$
By \cite[Ex 11.26]{rockwets}, we have
$$-e_{r+\delta(\lambda)}(F_\lambda)+\frac{r+\delta(\lambda)}{2}q=\left(F_\lambda+\frac{r+\delta(\lambda)}{2}q\right)^*((r+\delta(\lambda))\cdot),$$
where $f^*(x):= \sup_y\{\langle x, y \rangle - f(y)\}$ is the Fenchel conjugate as defined in \cite{howto}. This is an affine function composed with a convex function (as conjugate functions are convex), and as such it is convex. Notice that shifting the argument of $q$ by $\bar{x}$ only results in the addition of a linear term, as
$$q(x-\bar{x})=q(x)+2\langle x,\bar{x}\rangle+q(\bar{x})$$
where $q(\bar{x})$ is constant and $2\langle x,\bar{x}\rangle$ is linear. Hence, $PA_r+q(\cdot-\bar{x})$ is convex.\qed
\begin{thm}\textbf{[Stability of $PA_r$]}
\label{thm2}
For $i\in\{1,2,\ldots,m\}$, let $f_i:\mathbb{R}^n\rightarrow\mathbb{R}\cup\{\infty\}$ be proper, lsc, and prox-bounded with threshold $r_i.$ Let $\bar{r}>\max\limits_i\{r_i\}$ and $\bar{r} > \rho'$ from Theorem \ref{thm1} Condition 3. Suppose that for all $i,$ $P_{\bar{r}}f_i$ is single-valued and Lipschitz continuous (as is the case when $f_i$ is prox-regular). Then $PA_{\bar{r}}$ is well-defined and lower-$\mathcal{C}^2.$ If in addition
\begin{equation}
\label{eq1}
\Lip\left(\sum\limits_{i=1}^m\lambda_iP_{\bar{r}}f_i-I\right)\leq1,
\end{equation}
then for any $\bar{\lambda}$ such that $\delta(\bar{\lambda})>0$ we have
\begin{enumerate}
\item $PA_{\bar{r}}(\cdot,\bar{\lambda})\in\mathcal{C}^{1+}$ as a function of $x$
\item $PA_{\bar{r}}$ is locally Lipschitz continuous in $\lambda$ near $\bar{\lambda}$
\item $\nabla_xPA_{\bar{r}}$ is locally Lipschitz continuous in $\lambda$ near $\bar{\lambda}.$
\end{enumerate}
Finally, if $f_i+\frac{\bar{r}}{2}q$ is convex then $PA_{\bar{r}}(\cdot,e_i)=f_i(\cdot).$
\end{thm}
\textbf{Proof:} Let $F(x,\lambda)=-\sum\limits_{i=1}^m\lambda_ie_{\bar{r}}f_i(x).$ By Proposition \ref{prop1}, $PA_{\bar{r}}$ is well-defined and finite-valued. Since $P_{\bar{r}}f_i$ is single-valued for all $i,$ $P_{\bar{r}}F$ is single-valued as well. Since $f_i$ is proper, lsc, and prox-bounded for all $i,$ and $\bar{r}$ is greater than each threshold $r_i,$ Proposition \ref{cor1} gives us that $F$ is continuously para-prox-regular at $(\bar{x},\bar{\lambda})$ for $\bar{v}\in\partial_xF(\bar{x},\bar{\lambda}),$ and that $(0,y)\in\partial^\infty F(\bar{x},\bar{\lambda})\Rightarrow y=0.$ Since $P_{\bar{r}}f_i$ is single-valued and Lipschitz continuous for all $i,$ we have all the conditions of \cite[Prop 4.5]{proxave}, and therefore
\begin{enumerate}
    \item $(0,\lambda')\in D^*(\partial_xF)(\bar{x},\bar{\lambda}|\bar{v})(0)\Rightarrow\lambda'=0$
    \item $(x',\lambda')\in D^*(\partial_xF)(\bar{x},\bar{\lambda}|\bar{v})(v'), v'\neq0\Rightarrow\langle x',v'\rangle>-\rho|v'|^2$ for some $\rho>0$
    \item The mapping $\partial_xF(\bar{x},\cdot)$ has a continuous selection $g$ near $\bar{\lambda}.$
\end{enumerate} Hence the condition $\bar{r}>\max\{\rho,\rho',r\}$ of Theorem \ref{thm1} is satisfied (recall $r=\max_i\{r_i\}$). Therefore, all conditions of Theorem \ref{thm1} hold, and we may assume its result. Since $\delta\in\mathcal{C}^2,$ there exists $\bar{K}>0$ such that
$$|\delta(\lambda')-\delta(\lambda)|\leq\bar{K}|\lambda'-\lambda|$$
for all $\lambda',\lambda$ near $\bar{\lambda}.$ The rest of the proof is the same as that of \cite[Thm 4.6]{proxave}.\qed
\begin{cor}
\label{cor3}
For $i\in\{1,2,\ldots,m\}$, let $f_i:\mathbb{R}^n\rightarrow\mathbb{R}\cup\{\infty\}$ be proper and lsc such that for some $r>0,$ $f_i+\frac{r}{2}q$ is convex for all $i.$ Then $f_i$ is prox-regular and prox-bounded, and inequality (\ref{eq1}) holds. In particular, all the conditions of Theorem \ref{thm2} hold.
\end{cor}
\textbf{Proof:} Since $f_i+\frac{r}{2}q$ is convex for all $i,$ we have that $f_i$ is prox-bounded and lower-$\mathcal{C}^2$, and therefore prox-regular, for all $i.$ Since
$$P_1(f_i+\frac{r}{2}q)=P_{r+1}f_i,$$
by \cite[Prop 12.19]{rockwets} we have that $I-P_1(f_i+\frac{r}{2}q)$ is Lipschitz continuous with constant at most 1. Thus
$$\Lip\left\{\sum\limits_{i=1}^m\lambda_iP_{r+1}f_i-I\right\}=\Lip\left\{\sum\limits_{i=1}^m\lambda_i(I-P_{r+1}f_i)\right\}\leq\sum\limits_{i=1}^m\lambda_1=1.$$
This provides inequality (\ref{eq1}).\qed

\section{Example}
\label{example}
In 2010, Goebel, Hare, and Wang presented a study of the minimizers of the proximal average function for convex functions.   For convex functions $f_i$ recall that\\ $-e_r \left(-\sum_{i=1}^m \lambda_i e_r f_i\right)(x)$ defined the proximal average from \cite{howto}.  It was shown that
    \[\Phi(\lambda) := \argmin_x -e_r \left(-\sum_{i=1}^m \lambda_i e_r f_i\right)(x)\]
is single-valued and continuous, provided that all functions are bounded below and at least one function is essentially strictly convex \cite[Thm 3.8]{optval}. We next show that if $f_i$ are convex functions, then the minimizers of the NC-proximal average coincide exactly with the minimizers of the proximal average. In particular, in this case all results from \cite{optval} hold.
\begin{lem}\label{lem:argmin}For $i\in\{1,2,\ldots,m\}$ let $f_i:\mathbb{R}^n\rightarrow\mathbb{R}\cup\{\infty\}$ be proper, lsc, convex, and bounded below. Let $\lambda \in \Lambda$, then
    \[\argmin_x PA_r(x,\lambda) = \argmin_x \sum_{i=1}^m \lambda_i e_r f_i(x) = \argmin_x -e_r \left(-\sum_{i=1}^m \lambda_i e_r f_i\right)(x).\]
\end{lem}
\textbf{Proof:} The minimizers of $PA_r(\cdot,\lambda)$ coincide with the minimizers of its Moreau envelope $e_{r+\delta(\lambda)} PA_r(\cdot,\lambda)$.  By \cite[Ex 11.26(d)]{rockwets}, we have that $-e_{r+\delta(\lambda)} PA_r(x,\lambda) = \left(\sum_{i=1}^m -\lambda_i e_r f_i(x)\right)$, so the first equality holds.  The second equality appears in \cite[Lem 3.2]{optval}. \medskip\qed\\
If $f_i$ are non-convex, then the proximal average is undefined, and the results from \cite{optval} no longer apply.  In this case, the results of Theorem \ref{thm2} provide some small understanding of the continuity of the minimizers of the NC-proximal average, as follows.
\begin{cor}
\label{cor4}
Let the conditions of Theorem \ref{thm2} hold. Let $x_k\in\argmin\limits_xPA_r(x,\lambda_k).$ Suppose $\lambda_k\rightarrow\bar{\lambda}$ and $x_k\rightarrow\bar{x}.$ Then $\nabla PA_r(\bar{x},\bar{\lambda})=0.$
\end{cor}
\textbf{Proof:} By Theorem \ref{thm2}, $\nabla PA_r$ is Lipschitz continuous in $\lambda.$ Therefore, there exists $c>0$ such that for all $k,$
$$|\nabla PA_r(x_k,\lambda_k)-\nabla PA_r(x_k,\bar{\lambda})\leq c|\lambda_k-\bar{\lambda}|.$$
Since $x_k\in\argmin PA_r(x_k,\lambda_k),$ we know that $\nabla PA_r(x_k,\lambda_k)=0.$ So for all $k,$
$$|\nabla PA_r(x_k,\bar{\lambda})|\leq c|\lambda_k-\bar{\lambda}|.$$
Taking the limit as $k\rightarrow\infty,$ we find that $\nabla PA_r(\bar{x},\bar{\lambda})=0$.\medskip\qed\\
While Corollary \ref{cor4} gives us a way to identify the minimizers of $PA_r,$ it says nothing about the single-valuedness or the continuity of said minimizers. The example that follows illustrates that, in fact, the function of minimizers of the NC-proximal average may be multi-valued and discontinuous.\medskip

\begin{figure}[ht]
\begin{center}\includegraphics[scale=0.22]{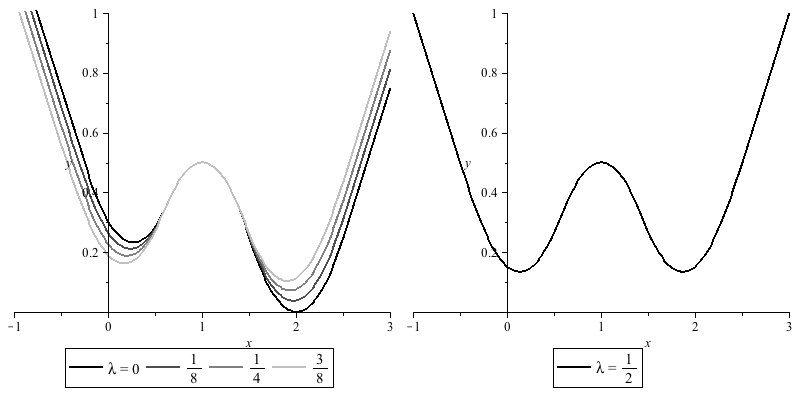}\end{center}
\caption{$G(\bar{x},\lambda)$}
\label{fig:2}
\end{figure}

Let $\epsilon=\frac{1}{2},$ and define the functions $g_0$ and $g_1$ via
\begin{align*}
 &g_0(x):=\max\{-x,-\frac{1}{2}(x-1)^2+\frac{1}{2},x-2+\epsilon\},\\
 &g_1(x):=\max\{-x+\epsilon,-\frac{1}{2}(x-1)^2+\frac{1}{2},x-2\}.
\end{align*}
Then $g_0$ and $g_1$ are proper, lsc, and bounded below.
\begin{figure}[h!]
\begin{center}\includegraphics[scale=0.5]{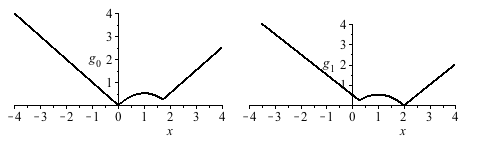}\end{center}
\caption{Functions $g_0$ and $g_1$ for $\epsilon=0.5.$}
\label{fig:1}
\end{figure}
Moreover, $g_i+\frac{1}{2}q$ is convex for $i\in\{1,2\}$. Let
$k=2-\sqrt{4-2\epsilon},~l=\sqrt{4-2\epsilon}$ and define
$$\begin{array}{l l l l}
\delta_0:=0 & \delta_1:=\epsilon & \epsilon_0:=\epsilon & \epsilon_1:=0\\
k_0:=0 & k_1:=k & l_0:=l & l_1:=2.
\end{array}$$
Consider $P_rg_i(\bar{x})=\argmin\limits_x\{g_i(x)+\frac{r}{2}|x-\bar{x}|^2\}.$ If $r > 1$, then we find that
$$P_r g_i(\bar{x})=\begin{cases}
\bar{x}+\frac{1}{r}, & \bar{x}<k_i-\frac{1}{r}\\
k_i, & \bar{x}\in[k_i-\frac{1}{r},k_i-\frac{k_i}{r}+\frac{1}{r}]\\
\frac{r\bar{x}-1}{r-1}, & \bar{x}\in(k_i-\frac{k_i}{r}+\frac{1}{r},l_i-\frac{l_i}{r}+\frac{1}{r})\\
l_i, & \bar{x}\in[l_i-\frac{l_i}{r}+\frac{1}{r},l_i+\frac{1}{r}]\\
\bar{x}-\frac{1}{r}, & \bar{x}>l_i+\frac{1}{r}.
\end{cases}$$
Evaluating the Moreau envelope and simplifying, we get
$$e_rg_i(\bar{x})=\begin{cases}
-\bar{x}-\frac{1}{2r}+\delta_i, & \bar{x}<k_i-\frac{1}{r}\\
\frac{r}{2}\bar{x}^2-rk_i\bar{x}+\frac{r-1}{2}k_i^2+k_i, & \bar{x}\in[k_i-\frac{1}{r},k_i-\frac{k_i}{r}+\frac{1}{r}]\\
-\frac{1}{2(r-1)}(r\bar{x}^2-2r\bar{x}+1), & \bar{x}\in(k_i-\frac{k_i}{r}+\frac{1}{r},l_i-\frac{l_i}{r}+\frac{1}{r})\\
\frac{r}{2}\bar{x}^2-rl_i\bar{x}+\frac{r-1}{2}l_i^2+l_i, & \bar{x}\in[l_i-\frac{l_i}{r}+\frac{1}{r},l_i+\frac{1}{r}]\\
\bar{x}-2-\frac{1}{2r}+\epsilon_i, & \bar{x}>l_i+\frac{1}{r}.
\end{cases}$$
Considering the specific example $r=2$, and applying $\epsilon=\frac{1}{2}$, we define the function $G(\bar{x},\lambda):=(\lambda e_2 g_0+(1-\lambda)e_2 g_1)(\bar{x}),$ which can be expanded to
$$G(\bar{x},\lambda)=\begin{cases}
-\bar{x}-\frac{\lambda}{2}+\frac{1}{4}, & x<-\frac{1}{2}\\
\lambda \bar{x}^2+(\lambda-1)\bar{x}-\frac{\lambda-1}{4}, & x\in[-\frac{1}{2},\frac{3-2\sqrt{3}}{2})\\
\bar{x}^2+(\lambda-1)(4-2\sqrt{3})\bar{x}-\frac{(\lambda-1)(11-6\sqrt{3})}{2}, & x\in[\frac{3-2\sqrt{3}}{2},\frac{1}{2}]\\
(1-2\lambda)\bar{x}^2+[-4+2\sqrt{3}+(6-2\sqrt{3})\lambda]\bar{x}+\frac{11-6\sqrt{3}}{2}-(6-3\sqrt{3})\lambda, & x\in(\frac{1}{2},\frac{3-\sqrt{3}}{2}]\\
-\bar{x}^2+2\bar{x}-\frac{1}{2}, & x\in(\frac{3-\sqrt{3}}{2},\frac{1+\sqrt{3}}{2})\\
(2\lambda-1)\bar{x}^2+[2-(2+2\sqrt{3})\lambda]\bar{x}-\frac{1}{2}+(2+\sqrt{3})\lambda, & x\in[\frac{1+\sqrt{3}}{2},\frac{3}{2})\\
\bar{x}^2-[4-(4-2\sqrt{3})\lambda]\bar{x}+4-\frac{5-2\sqrt{3}}{2}\lambda, & x\in[\frac{3}{2},\frac{1+2\sqrt{3}}{2}]\\
(1-\lambda)\bar{x}^2+(5\lambda-4)\bar{x}+4-\frac{23}{4}\lambda, & x\in(\frac{1+2\sqrt{3}}{2},\frac{5}{2}]\\
\bar{x}+\frac{\lambda}{2}-\frac{9}{4}, & x<\frac{5}{2}.
\end{cases}$$
By Lemma \ref{lem:argmin}, we know that
$$\argmin\limits_{\bar{x}}PA_r(\bar{x},\lambda)=\argmin\limits_{\bar{x}}G(\bar{x},\lambda).$$
Figure \ref{fig:2} displays graphs of $G$ for various values of $\lambda.$

Noting that $G\in\mathcal{C}^1,$ we find three critical points (where $\frac{\partial}{\partial x}G(x, \lambda)=0$):
\begin{enumerate}
\item $\bar{x}_1=(1-\lambda)(2-\sqrt{3})$ (leftmost local minimum argument),
\item $\bar{x}_2=1$ (local maximum argument),
\item $\bar{x}_3=2-(2-\sqrt{3})\lambda$ (rightmost local minimum argument).
\end{enumerate}
Observe that when $\lambda=\frac{1}{2}$ we have that $\bar{x}_1=\frac{2-\sqrt{3}}{2}$, $\bar{x}_3=\frac{2+\sqrt{3}}{2},$ and
$$G(\frac{2-\sqrt{3}}{2},\frac{1}{2})=\frac{2-\sqrt{3}}{2}=G(\frac{2+\sqrt{3}}{2},\frac{1}{2}).$$
This verifies that there are two minimizers when $\lambda=\frac{1}{2}.$ Finally, we note that
\begin{align*}
G(\bar{x}_1,\lambda)<G(\bar{x}_3,\lambda),~\lambda\in[0,\frac{1}{2})&&\mbox{ and }&&
G(\bar{x}_1,\lambda)>G(\bar{x}_3,\lambda),~\lambda\in(\frac{1}{2},1],
\end{align*}which proves the argmin is a singleton whenever $\lambda\neq\frac{1}{2}.$ Therefore, $\argmin PA_r$ is not a continuous function of $\lambda.$
\section{Conclusion}
We have seen that, using the Moreau envelope definition, the NC-proximal average can be generalized to accomodate any finite number of suitable functions. Under appropriate conditions, $PA_r$ is well-defined, lower-$\mathcal{C}^2,$ and locally Lipschitz continuous in $x$ and in $\lambda.$ These properties make $PA_r$ a useful function for researchers in the Optimization field.
\bibliographystyle{plain}
\bibliography{Bibliography}

\end{document}